# A Machine Learning Framework for Computing the Most Probable Paths of Stochastic Dynamical Systems

Yang Li[1,2,a], Jinqiao Duan[2,b] and Xianbin Liu[1,c]

[1]State Key Laboratory of Mechanics and Control of Mechanical Structures, College of Aerospace Engineering, Nanjing University of Aeronautics and Astronautics, 29 Yudao Street, Nanjing 210016, China

[2]Department of Applied Mathematics, College of Computing, Illinois Institute of Technology, Chicago, Illinois 60616, USA

[a]li_yang@nuaa.edu.cn

[b]duan@iit.edu

[c]Corresponding author: xbliu@nuaa.edu.cn

**Abstract** The emergence of transition phenomena between metastable states induced by noise plays a fundamental role in a broad range of nonlinear systems. The computation of the most probable paths is a key issue to understand the mechanism of transition behaviors. Shooting method is a common technique for this purpose to solve the Euler-Lagrange equation for the associated action functional, while losing its efficacy in high-dimensional systems. In the present work, we develop a machine learning framework to compute the most probable paths in the sense of Onsager-Machlup action functional theory. Specifically, we reformulate the boundary value problem of Hamiltonian system and design a neural network to remedy the shortcomings of shooting method. The successful applications of our algorithms to several prototypical examples demonstrate its efficacy and accuracy for stochastic systems with both (Gaussian) Brownian noise and (non-Gaussian) Lévy noise. This novel approach is effective in exploring the internal mechanisms of rare events triggered by random fluctuations in various scientific fields.

# 1. Introduction

The interaction between nonlinearity and randomness in dynamical systems may lead to emergence of novel behaviors, which has no analogue in the deterministic case. In particular, the last two decades have witnessed an increasing number of investigations in the analysis of transition phenomena induced by noise in various scientific fields such as biology [1–4], physics [5–7], chemistry [8,9] and engineering [10]. To explore the mechanism of transition between metastable states is a challenging and pivotal task in stochastic dynamical systems.

The Onsager-Machlup (OM) action functional [11] is a critical tool to study the transition of stochastic dynamical systems. The idea is to represent the probability of a single path by a tube around it with fixed thickness. Based on path integral formulation, the probability of this tube can be estimated by the Onsager-Machlup (OM) action functional. Thus the complicated computations of probability are transformed to the variational problem of the OM functional whose minimizer is called the most probable transition path.

In 1953, Onsager and Machlup [11] firstly derived the OM functional for diffusion processes with linear drift and constant diffusion coefficients. The extension to nonlinear systems was subsequently undertaken by Tisza and Manning [12]. Additionally, there was another approach to deduce the OM functional developed by Dürr and Bach [13] with the application of Girsanov transformation for measures induced by diffusion processes. Chao and Duan [14] generalized this method to solve the more complex cases for the stochastic dynamical systems under (non-Gaussian) Lévy noise as well as (Gaussian) Brownian noise. Tang et al. [15,16] further derived the OM functional from another aspect for the overdamped Langevin equation with multiplicative Gaussian noise and for general stochastic interpretation.

According to variational principle, the most probable path connecting two states satisfies either Euler-Lagrange equation or the corresponding Hamiltonian system for a given OM function. The shooting method [17] is common technique to deal with this two-point boundary value problem of a system ordinary differential equations. However, this shooting method is inefficient and even ill-posed in high-dimensional systems [18].

In this article, we will employ a neural network to bypass the drawbacks of the shooting method. As a powerful approach of machine learning, neural networks have been extensively applied to solve the ordinary [19] and partial [20–23] differential equations, and to learn the governing laws from

data [24,25].

More specifically, the most probable transition path connecting the points $\mathbf{x}(0) = \mathbf{x}_0$ and $\mathbf{x}(T) = \mathbf{x}_f$ can be transformed into the following boundary value problem

$$\dot{\mathbf{x}} = \frac{\partial H}{\partial \mathbf{p}}, \quad \mathbf{x}(0) = \mathbf{x}_0,$$

$$\dot{\mathbf{p}} = -\frac{\partial H}{\partial \mathbf{x}}, \quad \mathbf{p}(T) = \lambda,$$

where the meaning and derivation process of this equation will be explained in Section 2. However, the boundary value of momentum $\lambda$ is unknown. Thus one needs to infer a good prior for $\lambda$, which depends on the position $\mathbf{x}_f$: each end point $\mathbf{x}_f$ of a path has different momentum $\lambda$ as the best prior. Then a neural network is used to map out the function $\lambda(\mathbf{x}_f)$. To generate the training set for the neural network, we randomly generated many $\lambda$ and get its correspondent $\mathbf{x}_f$ by a iterative algorithm (Algorithm 1 in Section 3). Such a training set was used to train the neural network as a function $\lambda(\mathbf{x}_f)$. Therefore, instead of manually adjusting $\lambda$ to manage the path to hit $\mathbf{x}_f$ in shooting method, the learned neural network can bypass these fussy operations by automatically outputting the best prior $\lambda$.

This article is arranged as follows. In Section 2, we briefly introduce the Onsager-Machlup theory and reset the boundary value problem for the corresponding Hamiltonian system. Then we reformulate the boundary conditions of the Hamiltonian system and design numerical schemes to compute the most probable transition path in Section 3. We test our method by numerical experiments in Section 4, and finally conclude with Discussion in Section 5.

## 2. Onsager-Machlup theory

Consider the following $n$-dimensional stochastic dynamical system

$$d\mathbf{x}(t) = \mathbf{f}(\mathbf{x}(t))dt + C(\mathbf{x}(t))d\mathbf{B}_t + d\mathbf{L}_t, \tag{1}$$

where $\mathbf{B}_t = [B_{1,t}, \cdots, B_{n,t}]^T$ is $n$-dimensional Brownian motion and $\mathbf{L}_t = [L_{1,t}, \cdots, L_{n,t}]^T$ is non-Gaussian Lévy motion (see Appendix) with the jump measure $\nu$ satisfying $\int_{|\mathbf{y}|<1} |\mathbf{y}| \nu(d\mathbf{y}) < \infty$. The vector $\mathbf{f}(\mathbf{x}) = [f_1(\mathbf{x}), \cdots, f_n(\mathbf{x})]^T$ is the drift coefficient and $CC^T$ denotes diffusion matrix.

For Eq. (1) with multiplicative Gaussian noise, an ambiguity in choosing the integration method leads to different stochastic interpretations and a general notation is the $\alpha$-interpretation [26]. In order to avoid confusion with the parameter $\alpha$ of $\alpha$-stable Lévy noise, we use $\kappa$-interpretation instead of $\alpha$. The values $\kappa=0$, $\kappa=1/2$ and $\kappa=1$ correspond to Ito's, Stratonovich's, and anti-Ito's, respectively. By modifying the drift term, the stochastic differential equation (1) can be transformed into the unified Stratonovich form

$$d\mathbf{x}(t) = \mathbf{b}(\mathbf{x}(t))dt + C(\mathbf{x}(t))d\mathbf{B}_t + d\mathbf{L}_t, \tag{2}$$

where $\mathbf{b}(\mathbf{x}) = \mathbf{f}(\mathbf{x}) + (\kappa - 1/2)\Delta \mathbf{f}$ denotes the modified drift coefficient with $\Delta f_i = \sum_j \sum_k C_{kj} \frac{\partial C_{ij}}{\partial x_k}$.

Under random disturbances, transitions between metastable states may occur. During this process, what we are mostly interested in is to determine the most probable transition path. Since the probability of a single path is zero, we turn to explore the probability that the stochastic trajectory passes through the neighborhood or tube of a certain path. Under the condition of the given thickness of the tube, the probability of trajectory staying in the tube actually describes the possibility that this specific path realizes. More precisely, consider a tube surrounding a reference path $\varphi(t)$, $t \in [0,T]$. If for $\varepsilon$ sufficiently small, the probability of the solution process $\mathbf{x}(t)$ lying in this tube can be estimated in the following form

$$\mathbb{P}_{\mathbf{x}_0}\left\{\sup_{0 \le t \le T}|\mathbf{x}(t) - \varphi(t)| \le \varepsilon\right\} \propto C(\varepsilon)\exp\left\{-\frac{1}{2}\int_0^T \mathrm{OM}(\dot{\varphi},\varphi)dt\right\}, \tag{3}$$

then the integrand $\mathrm{OM}(\dot{\varphi},\varphi)$ is called the Onsager-Machlup (OM) function, or Lagrangian as in classical mechanics. We also refer $\int_0^T \mathrm{OM}(\dot{\varphi},\varphi)dt$ as the Onsager-Machlup (OM) action functional.

Accordingly, it is seen that the probability of the tube around $\varphi(t)$ is exponentially dominated by its OM functional. Thus the global minimum of the OM functional corresponds to the path with largest probability, i.e., the most probable path. Therefore, the computation for the probability is transformed into the variational problem of the OM functional. Based on the classical results of analytical mechanics, the most probable path, connecting the points $\mathbf{x}_0$ and $\mathbf{x}_f$, satisfies Euler-Lagrange equation

$$\frac{d}{dt}\left(\frac{\partial \mathrm{OM}(\dot{\mathbf{x}},\mathbf{x})}{\partial \dot{\mathbf{x}}}\right) - \frac{\partial \mathrm{OM}(\dot{\mathbf{x}},\mathbf{x})}{\partial \mathbf{x}} = 0 \tag{4}$$

with the boundary conditions $\mathbf{x}(0) = \mathbf{x}_0$ and $\mathbf{x}(T) = \mathbf{x}_f$. Here we have replaced the symbol $\varphi$ as $\mathbf{x}$. Since this is a second-order differential equation, it is more convenient for numerical integration to transform it into an equivalent Hamiltonian system

$$\dot{\mathbf{x}} = \frac{\partial H}{\partial \mathbf{p}}, \dot{\mathbf{p}} = -\frac{\partial H}{\partial \mathbf{x}}, \quad (5)$$
$$\mathbf{x}(0) = \mathbf{x}_0, \mathbf{x}(T) = \mathbf{x}_f.$$

Here, $H(\mathbf{x},\mathbf{p})$ is the Legendre transform of OM function and $\mathbf{p} = \frac{\partial \text{OM}(\dot{\mathbf{x}},\mathbf{x})}{\partial \dot{\mathbf{x}}}$ is called the momentum. The projection $\mathbf{x}(t)$ of the solution $(\mathbf{x}(t),\mathbf{p}(t))$ to the coordinate space provides the most probable path.

According to Refs. [14,15,27], the OM function of Eq. (2) is given, up to an additive constant, by

$$\text{OM}(\dot{\mathbf{x}},\mathbf{x}) = (\dot{\mathbf{x}}-\mathbf{b})^T (CC^T)^{-1} (\dot{\mathbf{x}}-\mathbf{b}) + Tr\left[C\nabla(C^{-1}\mathbf{b})\right]$$
$$+ 2(\dot{\mathbf{x}}-\mathbf{b})^T (CC^T)^{-1} \int_{|\mathbf{y}|<1} \mathbf{y}\nu(\mathrm{d}\mathbf{y}), \quad (6)$$

where $Tr[A] = \sum_i A_{ii}$. For convenience, denote $\mathbf{d} \triangleq \int_{|\mathbf{y}|<1} \mathbf{y}\nu(\mathrm{d}\mathbf{y})$. The Hamiltonian is then calculated by Legendre transform as

$$H(\mathbf{x},\mathbf{p}) = \frac{1}{4}\mathbf{p}^T CC^T \mathbf{p} + (\mathbf{b}-\mathbf{d})\cdot\mathbf{p} - Tr\left[C\nabla(C^{-1}\mathbf{b})\right] + \mathbf{d}^T (CC^T)^{-1} \mathbf{d}. \quad (7)$$

Hence the Hamiltonian system has the following form

$$\dot{\mathbf{x}} = \mathbf{b} - \mathbf{d} + 1/2\, CC^T \mathbf{p},$$
$$\dot{\mathbf{p}} = -1/4\nabla(\mathbf{p}^T CC^T \mathbf{p}) - (\nabla\mathbf{b})^T \mathbf{p} + \nabla\left\{Tr\left[C\nabla(C^{-1}\mathbf{b})\right]\right\} - \nabla\left[\mathbf{d}^T (CC^T)^{-1} \mathbf{d}\right], \quad (8)$$
$$\mathbf{x}(0) = \mathbf{x}_0, \mathbf{x}(T) = \mathbf{x}_f.$$

For the case of additive Gaussian noise, it is reduced as

$$\dot{\mathbf{x}} = \mathbf{f} - \mathbf{d} + 1/2\, CC^T \mathbf{p},$$
$$\dot{\mathbf{p}} = -(\nabla\mathbf{f})^T \mathbf{p} + \nabla(\nabla\cdot\mathbf{f}), \quad (9)$$
$$\mathbf{x}(0) = \mathbf{x}_0, \mathbf{x}(T) = \mathbf{x}_f.$$

Let us make two remarks on the form of OM function in Eq. (6). On one hand, if we consider $\alpha$-stable Lévy motion, then the condition $\int_{|\mathbf{y}|<1} |\mathbf{y}|\nu(\mathrm{d}\mathbf{y}) < \infty$ requires $0 < \alpha < 1$. On the other hand, the OM function of diffusion process is recovered when we set $\nu = 0$ in Eq. (6). In other words, the effect of Lévy noise on the transition phenomenon is reflected in the third term of OM function in Eq. (6). Specifically, if the Lévy motion is symmetric, then $\mathbf{d} = 0$ and the results are consistent with the ones of diffusion process consequently.

Note that Eq. (9) is a two-point boundary value problem of ordinary differential equation. Generally, it can be solved by shooting method. That is, we can adjust the initial value of the momentum and integrate the equation numerically until the final point $\mathbf{x}(T)$ reaches $\mathbf{x}_f$. However, there are still two shortcomings in this approach practically. First, we usually choose the stable fixed point called a metastable state as the initial point $\mathbf{x}_0$ to consider its transitions. Due to the conjugate momentum equation, i.e., the second equation of (9), containing the term $-(\nabla \mathbf{f})^T \mathbf{p}$, a numerical integration forward in time would be numerically unstable or even ill-posed. Moreover, it is a challenge in high-dimensional systems since it is hard to decide which component of the initial momentum should be adjusted. We will deal with these two shortcomings in next section.

## 3. Numerical algorithms

### 3.1. Reformulation

In order to deal with the divergence problem of the momentum, we generalize the method in Ref. [18] for Freidlin-Wentzell large deviation events to the case of the most probable paths under OM function with non-Gaussian Lévy noise. Denote $S_T[\mathbf{x}(t)] = \int_0^T \mathrm{OM}(\dot{\mathbf{x}}, \mathbf{x}) \mathrm{d}t$ and

$$I(\mathbf{x}_f) = \inf_{\mathbf{x} \in \mathcal{C}_1} S_T[\mathbf{x}(t)], \tag{10}$$

where $\mathcal{C}_1 = \{\mathbf{x} \in C[0,T] \mid \mathbf{x}(0) = \mathbf{x}_0, \mathbf{x}(T) = \mathbf{x}_f\}$. In addition, define

$$I^*(\lambda) = \inf_{\mathbf{x} \in \mathcal{C}_0} \{S_T[\mathbf{x}(t)] - \lambda \cdot \mathbf{x}(T)\} \tag{11}$$

with $\mathcal{C}_0 = \{\mathbf{x} \in C[0,T] \mid \mathbf{x}(0) = \mathbf{x}_0\}$. Note that this minimization does not require the constraint in the final point. In other words, the function space $\mathcal{C}_0$ describes the set of continuous trajectories starting at $\mathbf{x}_0$ regardless of their final point. In fact, $I^*(\lambda)$ and $I(\mathbf{x}_f)$ are Fenchel-Legendre duals. This can be motivated by

$$I^*(\lambda) = \inf_{\mathbf{x} \in \mathcal{C}_0} \{S_T[\mathbf{x}(t)] - \lambda \cdot \mathbf{x}(T)\}$$
$$= \inf_{\mathbf{x}_f \in R^n} \inf_{\mathbf{x} \in \mathcal{C}_1} \{S_T[\mathbf{x}(t)] - \lambda \cdot \mathbf{x}(T)\}$$
$$= \inf_{\mathbf{x}_f \in R^n} \left\{ \inf_{\mathbf{x} \in \mathcal{C}_1} S_T[\mathbf{x}(t)] - \lambda \cdot \mathbf{x}_f \right\}$$
$$= \inf_{\mathbf{x}_f \in R^n} \{I(\mathbf{x}_f) - \lambda \cdot \mathbf{x}_f\}.$$

Effectively, the connection between $I^*(\lambda)$ and $I(\mathbf{x}_f)$ allows us to deal with the minimization problem (11) instead of (10). Thus the variational results of (11) leads to the following Hamiltonian system

$$\dot{\mathbf{x}} = \frac{\partial H}{\partial \mathbf{p}}, \quad \mathbf{x}(0) = \mathbf{x}_0,$$
$$\dot{\mathbf{p}} = -\frac{\partial H}{\partial \mathbf{x}}, \quad \mathbf{p}(T) = \lambda. \tag{12}$$

It is seen from the expression of Eq. (5) that the boundary conditions are constrained on initial and final coordinates but absent on the momentum. By contrast, they are transformed into one boundary condition on coordinate and another on momentum in Eq. (12). The advantage of this operation is evident. The computation of Eq. (12) is to integrate $\mathbf{x}(t)$ forward in time and then to integrate the momentum backward in time. Thus we succeed in resolving the first shortcoming of the shooting method mentioned previously as the two directions of integration are both convergent.

The algorithm of this method is concluded in the following form:

**Algorithm 1:**

Step 1. Given a value of $\lambda$ and an initial guess trajectory $\mathbf{x}^k(t)$ ($k=1$ initially);

Step 2. Integrating the equation $\dot{\mathbf{p}} = -\frac{\partial H}{\partial \mathbf{x}}(\mathbf{x}^k, \mathbf{p})$ from $\mathbf{p}(T) = \lambda$ backward in time to obtain $\mathbf{p}^k$;

Step 3. Integrating the equation $\dot{\mathbf{x}} = \frac{\partial H}{\partial \mathbf{p}}(\mathbf{x}, \mathbf{p}^k)$ from $\mathbf{x}(0) = \mathbf{x}_0$ forward in time to obtain $\mathbf{x}^{k+1}$;

Step 4. Iterating Step 2-3 until convergence.

**3.2. Neural network**

Although the divergence problem of the momentum component is overcome by the proposed algorithm in the previous subsection, how to determine the value of $\lambda$ according to the final boundary condition of $\mathbf{x}$ still remains unsolved. This procedure can be accomplished by a deep neural network.

The architecture of this neural network is illustrated in Fig. 1. According to our purpose, the input and output states are fixed as $\mathbf{x}_f$ and $\lambda$ respectively, which have $n$ components each. Assume that there are $L$ hidden layers between input and output with $n_l$ neurons in $l$-th layer, $l = 1, 2, ..., L$. Accordingly, a complicated nonlinear function can be approximated by a composition of simpler functions

$$\lambda(\mathbf{x}_f) = g^{(L+1)} \circ g^{(L)} \circ \cdots \circ g^{(1)}(\mathbf{x}_f), \tag{13}$$

where each layer $g^{(l)}(\cdot)$ is defined as

$$g^{(l)}(\mathbf{a}^{(l-1)}) = \sigma^{(l)}(W^{(l)}\mathbf{a}^{(l-1)} + \mathbf{b}^{(l)}), l = 1, 2, ..., L+1. \tag{14}$$

Here we indicate $\mathbf{a}^{(0)} = \mathbf{x}_f$ and $\mathbf{a}^{(L+1)} = \lambda$ for convenience. Moreover, $W^{(l)}$ and $\mathbf{b}^{(l)}$ are called the weight matrices and bias vectors, respectively. And $\sigma^{(l)}$ represents a nonlinear activation function applied component-wise to its argument, of which popular choices include sigmoid, tanh, ReLu and other similar functions. In what follows, we use the ReLu function $\text{ReLu}(x) = \max\{0, x\}$ for all the hidden layers. The output layer is chosen as the identity function in order to cover the space of $\lambda$.

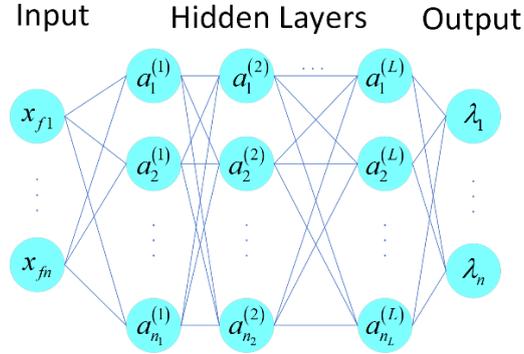

Fig. 1. Architecture of neural network with $L$ hidden layers.

Let $\theta$ denote the collection of the parameters of the neural network, i.e., $\theta = \{W^{(l)}, \mathbf{b}^{(l)} : l = 1, 2, ..., L+1\}$. The neural network is trained by optimizing over the parameters $\theta$ to best approximate $\lambda(\mathbf{x}_f)$ by its output states. Specifically, assume that we have the data set $\mathcal{D} = \{(\mathbf{x}_f^m, \lambda^m) : m = 1, 2, ..., M\}$ and introduce the loss function as the $L_2$-distance between the prediction of the network and the targets. Then the neural network is trained by searching the optimal parameters to solve the nonlinear regression problem of the cost function

$$J(\theta) = \frac{1}{2M} \sum_{m=1}^{M} \left[ \lambda^{NN}\left(\mathbf{x}_f^m; \theta\right) - \lambda^m \right]^2. \tag{15}$$

Generally, the minimizing parameters are iterated by the gradient descent method

$$\theta_{k+1} = \theta_k - \eta \left.\frac{\partial J(\theta)}{\partial \theta}\right|_{\theta=\theta_k}, \tag{16}$$

where the learning rate $\eta$ is a small number.

The training procedure of the neural network is summarized in the following form:

**Algorithm 2:**

Step 1. Given the data set and a group of initial guess parameters $\theta_0$;

Step 2. Using forward propagation to compute the output of the network and cost function;

Step 3. Using backward propagation to compute the gradients of the cost function to the parameters;

Step 4. Updating the parameters by Eq. (16);

Step 5. Iterating until convergence.

### 3.3. Algorithm

As a summary, the framework of our algorithm is listed as follows

**Algorithm 3:**

Step 1. Generating data. Specifically, we select $M$ points randomly in a suitable domain of the space of $\lambda$. With the initial point fixed as a stable state of the system, we compute the most probable path $\mathbf{x}(t)$ for every point of $\lambda$ in terms of Algorithm 1;

Step 2. Training neural network. Based on Algorithm 2, the data generated in step 1 are used to train the neural network with the input $\mathbf{x}(T)$ and output $\lambda$;

Step 3. Testing. Given a certain final point $\mathbf{x}_{test}$, we utilize the trained neural network to calculate the corresponding output $\lambda_{test}$;

Step 4. Computing the most probable path. We use the Algorithm 1 again to compute the most probable path for $\lambda_{test}$ and compare its end point with $\mathbf{x}_{test}$.

## 4. Numerical experiments

With the algorithm designed previously, we now present several prototypical examples to demonstrate our method for computing the most probable paths.

**Example 1** Consider a one-dimensional stochastic energy balance model [5] describing the climate change of earth

$$dx = -U'(x)dt + dB_t + dL_t \tag{17}$$

with the potential function

$$U(x) = \frac{1}{Ch}\left(-\frac{1}{4}S_0\left(0.5x + 2\ln\left(\cosh\frac{x-265}{10}\right)\right) + \frac{1}{5}\gamma\theta x^5\right), \tag{18}$$

where $x$ represents the global mean surface temperature. The Brownian motion $B_t$ and $\alpha$-stable Lévy motion $L_t$ are independent stochastic processes. In fact, the global energy change can be regarded as the difference between the incoming solar radiative energy and the outgoing radiative energy, $-U'(x) = 1/Ch(E_{in} - E_{out})$. Herein, the incoming energy $E_{in} = 1/4(1-\alpha(x))S_0$ with the planetary albedo $\alpha(x) = 0.5 - 0.2\tanh((x-265)/10)$ represents the total amount of solar radiation absorbed by the earth after the surface reflection. By considering the earth as a blackbody radiator with an effective surface temperature $x$, the outgoing energy is provided by Stefan-Boltzmann law $\theta x^4$. During this process, the greenhouse effect contributes to the global average temperature rising. Thus the outgoing energy is expressed as $E_{out} = \gamma\theta x^4$ with the greenhouse factor $\gamma \in [0,1]$. Note that the greenhouse effect increases as $\gamma$ decreases.

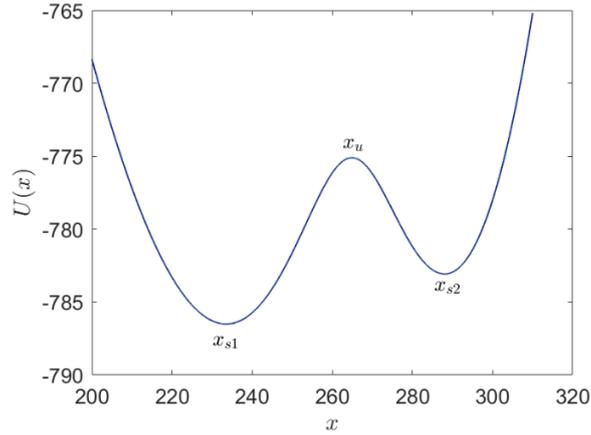

Fig. 2. The potential function of the energy balance model. The two minima correspond to the colder glacial state $x_{s1} = 233.52\text{K}$ (-39.63°C) and warmer interstadial state $x_{s2} = 288.03\text{K}$ (14.88°C).

In Eq. (18), the heat capacity $Ch = 46.8\text{Wyrm}^{-2}$ is defined as the amount of heat that must be added to the object in order to raise its temperature. Other parameters are chosen as the solar constant $S_0 = 1368 \text{Wm}^{-2}$ and the Stefan constant $\theta = 5.67 \times 10^{-8} \text{Wm}^{-2}\text{K}^{-4}$. The parameter $\gamma = 0.61$ corresponds to the bistable case with the colder glacial state $x_{s1} = 233.52\text{K}$ (-39.63°C) and warmer interstadial state $x_{s2} = 288.03\text{K}$ (14.88°C), separated by one unstable state $x_u$. The potential function of the system is shown in Fig. 2.

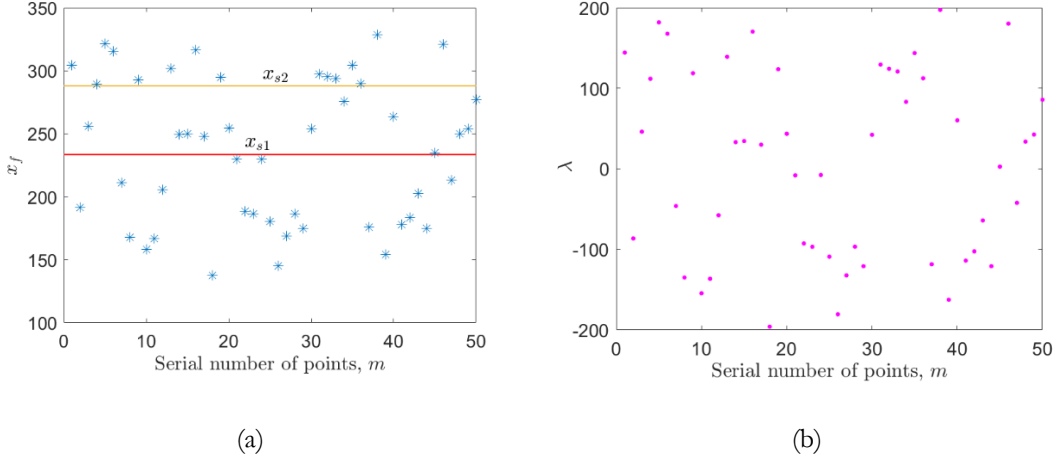

(a)            (b)

Fig. 3. (a) $M = 50$ integrated end points of $x_f$. (b) $M = 50$ points of randomly sampled $\lambda$.

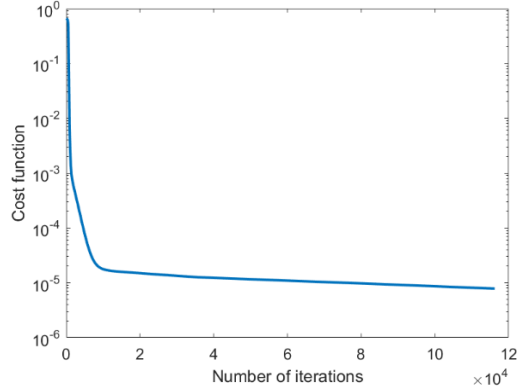

Fig. 4. The values of cost function of neural network during the training process in energy balance model.

Under random disturbances, we investigate the transition problem between the two stable states by applying our method. In this example, it is no need to distinguish different stochastic interpretations due to the additive Gaussian noise. In the Hamiltonian of Eq. (7), the drift term $b(x) = -U'(x)$, $C = 1$ and $d = \dfrac{\alpha\beta}{\Gamma(2-\alpha)\cos(\pi\alpha/2)}$ with the stability parameter $\alpha \in (0,1)$ and the skewness parameter $\beta \in [-1, 1]$. In what follows, we consider the case $\alpha = 0.5$ and $\beta = 0$. In addition, we fix the time interval length as $T = 1$ and choose $M = 50$ random points subject to uniform distribution in the domain $[-200, 200]$ of the space of $\lambda$. We use the neural network with 3 hidden layers and 20 neurons per layer. The learning rate is fixed as $\eta = 0.01$.

First, the Algorithm 1 is implemented within which the second-order Runge-Kutta method is adopted to integrate the Hamiltonian system from $x_0 = x_{s1}$ to obtain $M$ final points $x_f$, as shown in Fig. 3. In order to guarantee the convergence of the neural network, we rescale the input and output states as $y_f = \frac{2}{x_{s2} - x_{s1}} \left( x_f - \frac{x_{s1} + x_{s2}}{2} \right)$ and $\zeta = \lambda/200$. Second, the values of cost function are recorded during the training process and plotted in Fig. 4. It is seen that the network converges well and the error drops to the magnitude of $10^{-6}$.

As a result, the trained neural network is used to output the corresponding value of rescaled $\lambda$ when we input the rescaled $x_{s2}$. Subsequently we employ this value to compute the most probable transition path on the interval $[0, 1]$ by Algorithm 1, denoted as red dotted line in Fig. 5. It is found that the path reaches 287.73K at the end time which is sufficiently close to $x_{s2} = 288.03$K. Meanwhile, we also plot several paths iterated by Algorithm 1 with randomly sampled $\lambda$ in Fig. 5. It is observed that these paths are hard to reach our target $x_{s2}$ directly. This demonstrates the efficiency and advantage of the neural network to compute the most probable paths.

The results can be verified by two methods. On one hand, we can integrate the Hamiltonian system directly from the initial point $(x(0), p(0))$ where $p(0)$ is extracted from the computed path by our method. On the other hand, the most probable path can be directly obtained by minimizing the OM action functional $S_T[x(t)] = \int_0^T \text{OM}(\dot{x}, x) \text{d}t$. Numerically, we can divide the interval $[0, T]$ into $L$ subintervals with $0 = t_0 < t_1 < \cdots < t_L = T$. A path $x(t)$ is approximated by its values $\Phi_l$, at $t = t_l$ for $l = 0, 1, \ldots, L$. Then the action $S$ is approximated as

$$S_{t_0, \ldots, t_L}[\Phi_1, \ldots, \Phi_{L-1}] = \sum_{l=1}^{L} \text{OM}\left( \frac{\Phi_l - \Phi_{l-1}}{\Delta t_l}, \Phi_{l-1/2} \right) \Delta t_l, \tag{19}$$

where $\Delta t_l = t_l - t_{l-1}$ and $\Phi_{l-1/2} = (\Phi_l + \Phi_{l-1})/2$. Consequently, gradient descent method can be applied to minimize this objective function. Once the most probable path is determined, its conjugate momentum can be approximated by $p = \partial OM/\partial \dot{x}$. It is seen from Fig. 6 that the results of not only the most probable path but also the conjugate momentum agree well for the three methods.

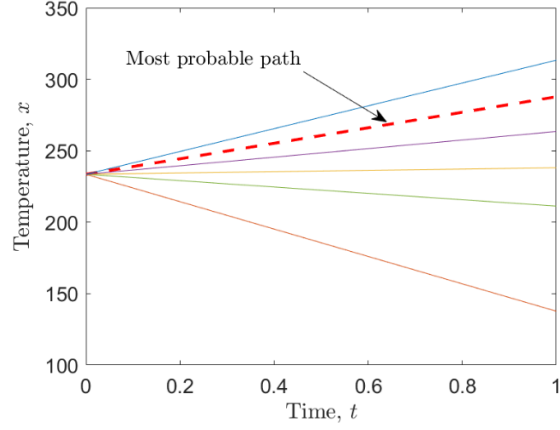

Fig. 5. Comparison of most probable paths between with using neural network and with randomly sampled $\lambda$. Red dotted curve denotes the most probable path with using neural network and other curves indicate the paths with random $\lambda$.

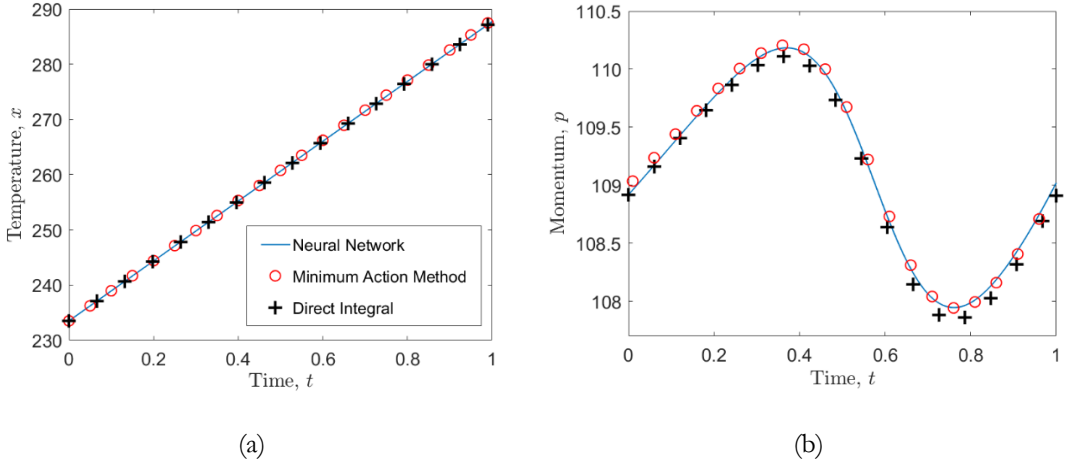

(a)          (b)

Fig. 6. (a) The most probable transition paths and (b) conjugate momenta, computed by neural network and Algorithm 1 (denoted as blue line), directly integrated by Hamiltonian system from initial point $(x(0), p(0))$ (indicated as black plus sign) and computed by minimum action method (shown as red circle).

**Example 2** Since the most important goal of our method is to deal with high-dimensional systems, the three-dimensional Lorenz system excited by both Gaussian Brownian noise and non-Gaussian Lévy noise is chosen as another example

$$d\mathbf{x}(t) = \mathbf{f}(\mathbf{x}(t))dt + C(\mathbf{x}(t))d\mathbf{B}_t + d\mathbf{L}_t,$$
$$\mathbf{f}(\mathbf{x}) = \left[\sigma(-x_1 + x_2), \rho x_1 - x_2 - x_1 x_3, -\gamma x_3 + x_1 x_2\right]^T, \quad (20)$$
$$C(\mathbf{x}) = \begin{bmatrix} 1 & 0 & 0 \\ 0 & 1 & 0 \\ 0 & 0 & \sqrt{1+\mu x_3^2} \end{bmatrix},$$

where the parameters in the model are fixed as $\sigma = 1$, $\gamma = 8/3$ and $\rho = 0.5$. $L_{1,t}$, $L_{2,t}$ and $L_{3,t}$ are independent $\alpha$-stable Lévy motions with $\alpha_1 = \alpha_2 = \alpha_3 = 0.5$. Since the diffusion matrix depends on state variable, the value of $\kappa$ controls the choice of stochastic interpretations, and the modified drift term is $\mathbf{b}(\mathbf{x}) = \mathbf{f}(\mathbf{x}) + (\kappa - 1/2)\Delta\mathbf{f}$ with $\Delta\mathbf{f} = [0, 0, \mu x_3]^T$.

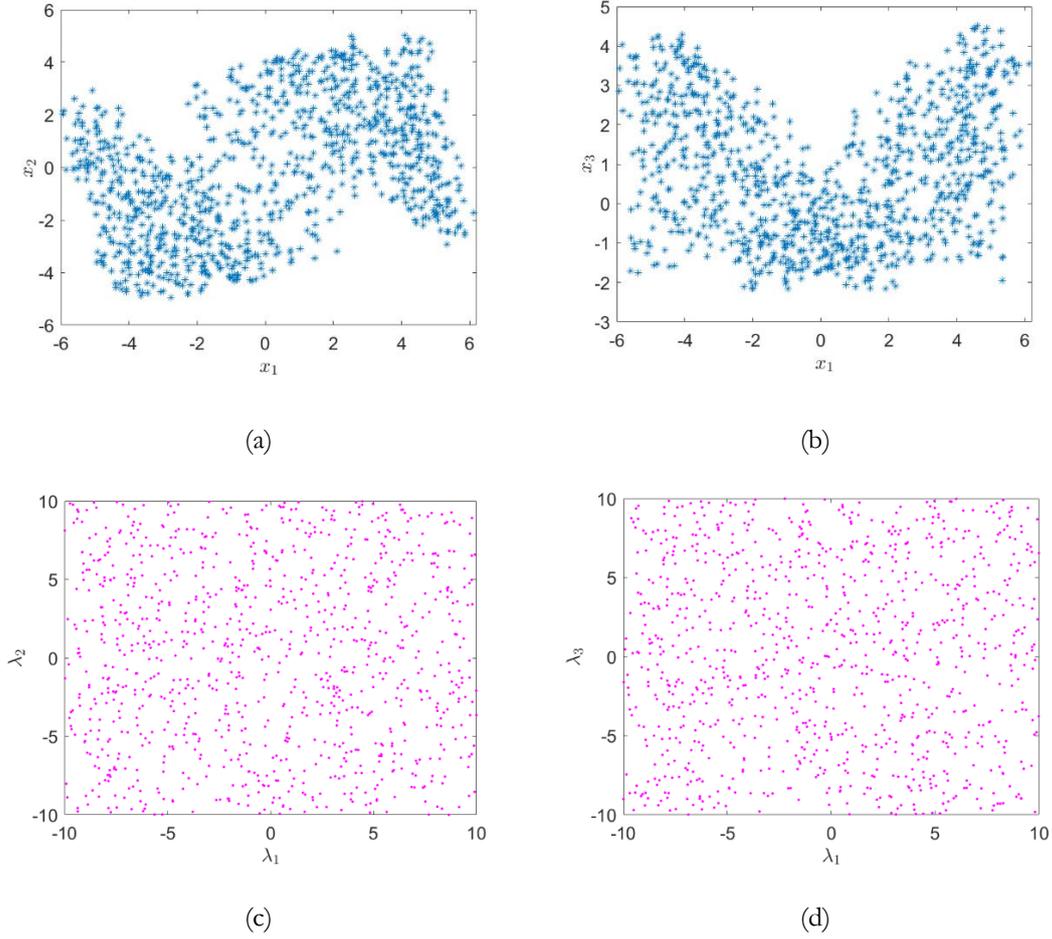

Fig. 7. $M = 1000$ randomly sampled $\lambda$ and corresponding integrated end points $\mathbf{x}_f$ for Case 1. (a) $x_{f2}$ versus $x_{f1}$. (b) $x_{f3}$ versus $x_{f1}$. (a) $\lambda_2$ versus $\lambda_1$. (a) $\lambda_3$ versus $\lambda_1$.

Note that the origin is unique stable fixed point. In what follows, we implement our approach to study the transition problem from origin to the point $(1,1,1)$. The time interval length is still fixed as $T=1$ and $M=1000$ random points are selected to cover the region $[-10,10] \times [-10,10] \times [-10,10]$ in the space of $\lambda$. We still use the neural network with 3 hidden layers and 20 neurons per layer with the learning rate $\eta = 0.01$. In what follows, we will consider four different cases.

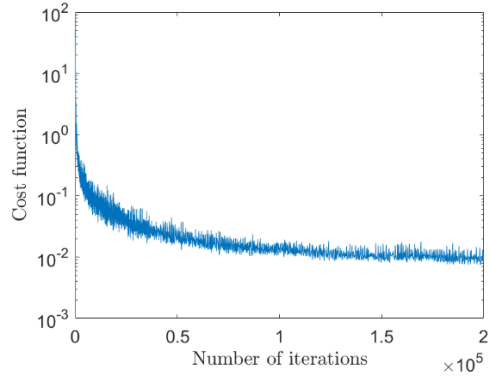

Fig. 8. The values of cost function of neural network during the training process for Case 1 in Lorenz model.

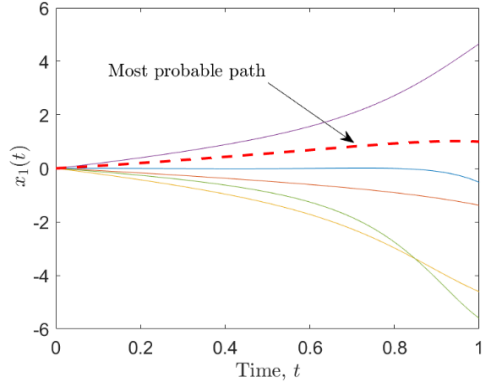

Fig. 9. Comparison of first component $x_1(t)$ of the most probable paths between with using neural network and with randomly sampled $\lambda$ for Case 1. Red dotted curve denotes the most probable path with using neural network and other curves indicate the paths with randomly sampled $\lambda$.

**Case 1**: $\kappa = 1/2$, $\mu = 0$, $\beta_1 = \beta_2 = \beta_3 = 0$.

This case is a relatively simple one with additive Gaussian noise since $\mu = 0$. As shown in Fig. 7, we first integrate the Hamiltonian system to obtain $M$ points $\mathbf{x}_f$ by performing Algorithm 1. Then the results are used to train the neural network. During the training process, the values of the cost function decrease to the magnitude of $10^{-2}$ and are plotted in Fig. 8. After the input of the point $(1,1,1)$ to the trained network and the implementation of Algorithm 1 again, the most probable path between the origin and $(1,1,1)$ is evaluated and its first component is illustrated in Fig. 9. The comparison to some integrated paths with randomly sampled $\lambda$ shows the efficiency of neural network in finding most probable path to reach a specific point.

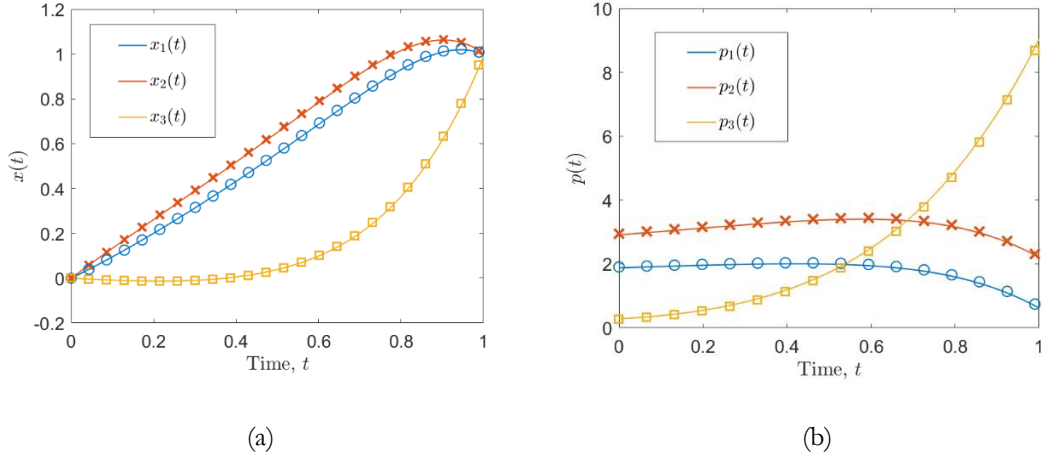

(a)        (b)

Fig. 10. Comparison of (a) the most probable transition path and (b) conjugate momentum between with neural network (denoted as lines) and with the minimum action method (indicated as circle, cross and square sign) for Case 1.

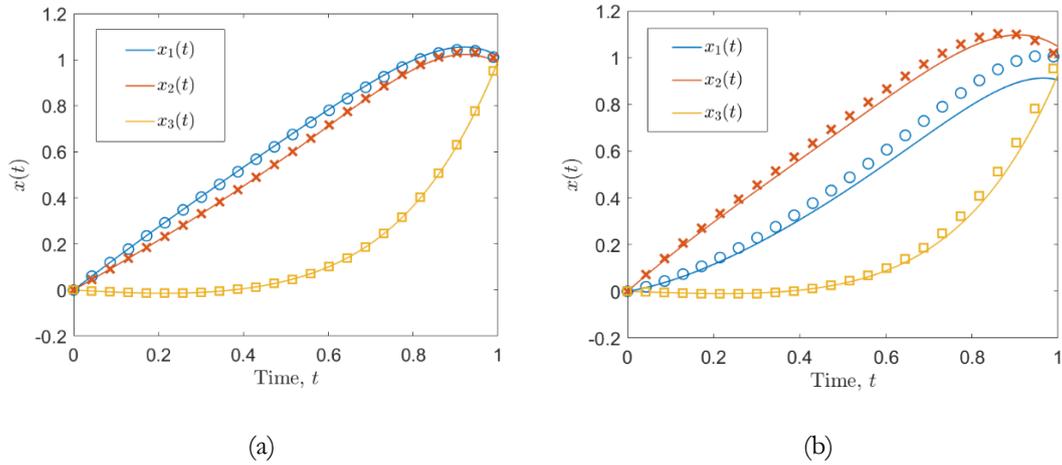

(a)        (b)

Fig. 11. Comparison of the most probable transition paths between with neural network (denoted as lines) and with the minimum action method (indicated as circle, cross and square sign) for Case 2 with (a) $\beta_1 = -0.5$ and (b) $\beta_1 = 0.5$.

The most probable transition path and its conjugate momentum are computed and shown in Fig. 10. Meanwhile, the minimum action method and gradient descent method are used to test our results.

It is found that the two results agree perfectly well, which implies that our algorithm is effective for additive Gaussian noise.

**Case 2**: $\kappa = 1/2$, $\mu = 0$, $\beta_1 = -0.5, 0.5$, $\beta_2 = \beta_3 = 0$.

In fact, Case 1 $\beta_1 = \beta_2 = \beta_3 = 0$ corresponds to the case that only Gaussian noise takes effect to drive the particle to transit. If we change $\beta_1$ to other values, then the Lévy noise of first direction starts to be effective. For instance, the results of $\beta_1 = -0.5$ and $\beta_1 = 0.5$ are exhibited in Fig. 11. It is found that their errors are slightly larger than the previous situation but within acceptable range. This implies that our method is more suitable to Gaussian noise than non-Gaussian noise. We can reduce the errors through the operations such as increasing the amount of data and the neurons in the network and decreasing the learning rate at the cost of computation time. In addition, it is seen that the result of $x_1(t)$ is tremendously affected than the other two components if we compare Figs. 10(a) and 11. Actually, it can be understood if we notice that the term $-\mathbf{d}$ which is dominated by $\beta$ emerges in the vector field of Eq. (9).

**Case 3**: $\kappa = 1/2$, $\mu = 1, 2$, $\beta_1 = \beta_2 = \beta_3 = 0$.

In view of the fact that $\mu = 0$ corresponds to additive Gaussian noise, we can regard the parameter $\mu$ as the deviation degree from additive to multiplicative noise. In order to reveal the impact of multiplicative noise to our numerical method, we compute the most probable paths for $\mu = 1, 2$ with neural network and with minimum action method as comparison and plot the results in Fig. 12. Combined with the results for $\mu = 0$ in Fig. 10(a), it is seen that this method will be less accurate with increasing $\mu$. For larger $\mu$, our method ceases to be effective since the Algorithm 1 may not converge for some region of $\lambda$. Therefore, our method is more effective for additive or weak multiplicative noise than strong multiplicative noise.

**Case 4**: $\kappa = 0, 1$, $\mu = 1$, $\beta_1 = \beta_2 = \beta_3 = 0$.

In order to reveal the impacts of different stochastic interpretations to our algorithm, we compute the most probable paths for $\kappa = 0, 1/2, 1$ corresponding to Ito's, Stratonovich's, and anti-Ito's sense. Fig. 13 shows the results for $\kappa = 0, 1$ with neural network and with minimum action method. Combined with $\kappa = 1/2$ in Fig. 12(a), it is found that the result for $\kappa = 1/2$ is better than $\kappa = 1$, and than $\kappa = 0$. Consequently, we can conclude that this method is most effective for Stratonovich's, then anti-Ito's, and then Ito's.

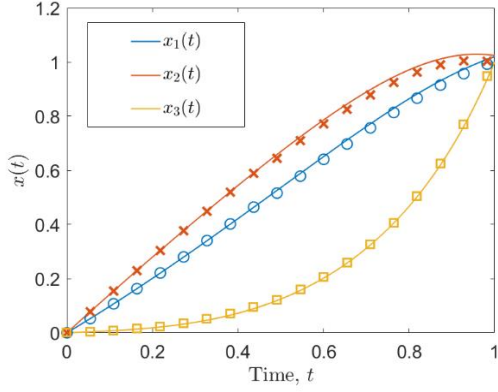 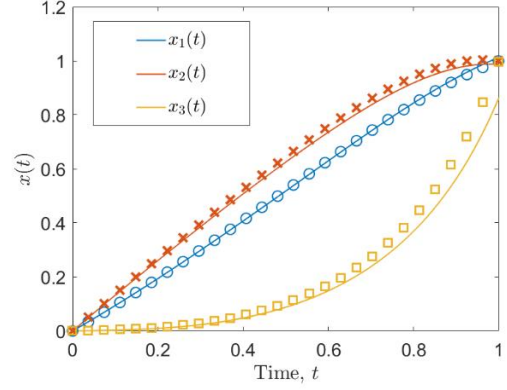

(a)                  (b)

Fig. 12. Comparison of the most probable transition paths between with neural network (denoted as lines) and with minimum action method (indicated as circle, cross and square sign) for Case 3 with (a) $\mu = 1$ and (b) $\mu = 2$.

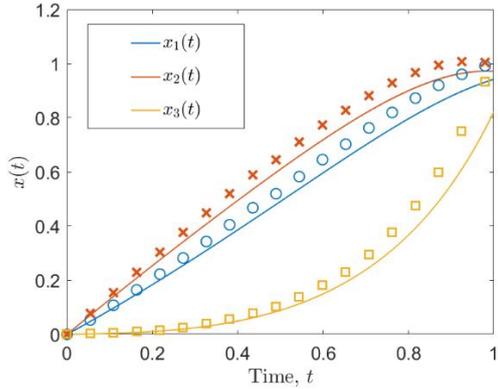 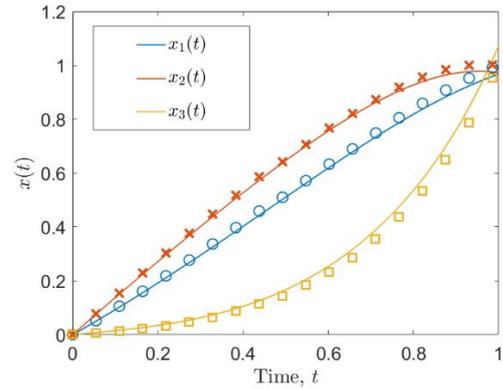

(a)                  (b)

Fig. 13. Comparison of the most probable transition paths between with neural network (denoted as lines) and with minimum action method (indicated as circle, cross and square sign) for Case 3 with (a) $\kappa = 0$ and (b) $\kappa = 1$.

## 5. Discussion

In this paper, we have devised a novel approach to compute the most probable transition paths of stochastic dynamical systems. Specifically, we briefly reviewed the Onsager-Machlup theory and derived the boundary value problem of the auxiliary Hamiltonian system. Then we reformulated its boundary condition and proposed the Algorithm 1 to compensate the divergence problem of the momentum direction in shooting method. Furthermore, we designed a neural network to automatically transform final boundary of the coordinate into the boundary condition of the momentum. The complete algorithm was concluded in Algorithm 3. The successful applications of our method to two prototypical examples confirmed its effectiveness for the systems with or without Lévy noise, additive or weak multiplicative Gaussian noise with Ito's, Stratonovich's, and anti-Ito's sense, and systems of various dimensions. Results show that our method is more suitable to Gaussian noise than non-Gaussian noise and more suitable to additive noise than multiplicative noise. Moreover, this algorithm is confirmed to be most effective for Stratonovich's, then anti-Ito's, and then Ito's.

Remark that our method can be generalized to other similar problems. For example, it can still work out in computing the most probable paths for Freidlin-Wentzell [28] large deviation theory in addition to Onsager-Machlup theory. More broadly, it can be used for other fields such as the optimal control problems.

According to our computed results, it is sufficiently accurate with 20 randomly generated points in one-dimensional system and with 500 points in three-dimensional system. Like other machine learning problems, it is inevitable to encounter a challenge in high dimensional system, which is called curse of dimensionality. Namely, the trained data information must increase exponentially with the number of dimension to guarantee the accuracy of estimate parameters. Thus it requires a huge amount of memory and impractical processor time. In order to overcome this issue, we can randomly generate some $\lambda$ points and integrate their correspondent $\mathbf{x}_f$. Then we can compute the distances between these $\mathbf{x}_f$ and our target, and take several closest points to find their preimages $\lambda$. Thus new $\lambda$ points can be randomly sampled in this smaller domain of $\lambda$ space which will cover our target due to the continuous feature of the map $\lambda(\mathbf{x}_f)$. This operation will tremendously reduce our data information required. On the other hand, the mini-batch or stochastic gradient descent method can be applied to update the parameters instead of the gradient descent method in order to save computing time.

Finally, it is worth noting that there still exists a challenge in the application of this algorithm. The neural network can be trained more effectively and accurately if the data $\mathbf{x}_f$ of integral results are relatively uniformly distributed such as Fig. 7. However, it does not work well if these points focus on a few specific regions, since the weights of various domains in cost function differ too much. Therefore, the successful application of our algorithm requires some kind of uniformity of vector field of the system.

**Acknowledgements**

We would like to thank Jianyu Hu, Ying Chao and Xiaoli Chen for helpful discussions. This research was supported by the National Natural Science Foundation of China (No. 11772149 and 1177144), A Project Funded by the Priority Academic Program Development of Jiangsu Higher Education Institutions (PAPD), The Research Fund of State Key Laboratory of Mechanics and Control of Mechanical Structures (MCMS-I-19G01), and the China Scholarship Council (CSC No. 201906830018).

**Appendix. The $\alpha$-stable Lévy motions**

A scalar $\alpha$-stable Lévy process $L_t$ is a stochastic process with the following conditions:

(i) $L_0 = 0$, a.s.;

(ii) Independent increments: for any choice of $n \geq 1$ and $t_0 < t_1 < \cdots < t_{n-1} < t_n$, the random variables $L_{t_0}$, $L_{t_1} - L_{t_0}$, $L_{t_2} - L_{t_1}$, $\cdots$, $L_{t_n} - L_{t_{n-1}}$ are independent;

(iii) Stationary increments: $L_t - L_s \sim S_\alpha \left( (t-s)^{1/\alpha}, \beta, 0 \right)$;

(iv) Stochastically continuous sample paths: for every $s > 0$, $L_t \to L_s$ in probability, as $t \to s$.

The $\alpha$-stable Lévy motion is a special but most popular type of the Lévy process defined by the stable random variable with the distribution $S_\alpha(\delta, \beta, \lambda)$ [29–31]. Usually, $\alpha \in (0, 2]$ is called the stability parameter, $\delta \in (0, \infty)$ is the scaling parameter, $\beta \in [-1, 1]$ is the skewness parameter and $\lambda \in (-\infty, \infty)$ is the shift parameter. The Lévy jump measure $\nu(\mathrm{d}y)$ has the following form

$$\nu(\mathrm{d}y) = \begin{cases} \dfrac{k_\alpha(1+\beta)}{2|y|^{1+\alpha}} \mathrm{d}y, & y > 0, \\ \dfrac{k_\alpha(1-\beta)}{2|y|^{1+\alpha}} \mathrm{d}y, & y < 0, \end{cases}$$

where

$$k_\alpha = \begin{cases} \dfrac{\alpha(1-\alpha)}{\Gamma(2-\alpha)\cos(\pi\alpha/2)}, & \alpha \neq 1, \\ \dfrac{2}{\pi}, & \alpha = 1. \end{cases}$$

A stable random variable $X$ with $0 < \alpha < 2$ has the following "heavy tail" estimate:

$$\lim_{x \to \infty} y^\alpha \mathbb{P}(X > y) = C_\alpha \frac{1+\beta}{2} \delta^\alpha;$$

where $C_\alpha$ is a positive constant depending on $\alpha$. In other words, the tail estimate decays polynomially. The $\alpha$-stable Lévy motion has larger jumps with lower jump frequencies for smaller $\alpha$ ($0 < \alpha < 1$), while it has smaller jump sizes with higher jump frequencies for larger $\alpha$ ($1 < \alpha < 2$). The special case $\alpha = 2$ corresponds to (Gaussian) Brownian motion. For more information about Lévy process, refer to Refs. [32,33].

## Data Availability Statement

The data that support the findings of this study are openly available in GitHub [34].